\documentclass[12pt]{article}
\def\hind{\hangindent=2pc\hangafter=1}
\newfont{\smcaps}{cmcsc10 scaled\magstep1}

\newcommand{\MA}{{\rm ~MA\,}}

\newcommand{\AR}{{\rm ~AR\,}}

\newcommand{\Cov}{{\rm ~Cov\,}}

\newcommand{\B}{{\cal B}}
\newcommand{\tr}{{\rm ~tr\,}}
\font\tty=cmtt10 at 11truept
\usepackage{epsfig}
\raggedright
\begin{document}
\baselineskip=22pt

\title{Faster ARMA Maximum Likelihood Estimation}
\date{}
\author{{A.I. McLeod$^{\rm a\/}$ and Y. Zhang$^{\rm b\/}$}\\\\
$^{\rm a\/}$\it Department of Statistical and Actuarial Sciences,\\
\it The University of Western Ontario,\\
\it London, Ontario Canada N6A 5B7\\\\
$^{\rm b\/}$\it Department of Mathematics and Statistics,\\
\it Acadia University,\\
\it Wolfville, Nova Scotia, Canada B4P 2R6
}
\maketitle
\hrule
\bigskip
\noindent {\bf Preprint: }\ A.I. McLeod andY. Zhang (2008), Faster ARMA maximum likelihood estimation, {\it Computational Statistics \& Data Analysis}, 52-4, 2166-2176.  doi: 10.1016/j.csda.2007.07.020

\newpage
{\bf Abstract\/}

\bigskip
\noindent
A new likelihood based AR approximation is given for ARMA models.
The usual algorithms for the computation of the likelihood of an ARMA model
require $O(n)$ flops per function evaluation.
Using our new approximation, an algorithm is developed which
requires only $O(1)$ flops in repeated likelihood evaluations.
In most cases, the new algorithm gives results identical to or very close to the exact
maximum likelihood estimate (MLE).
This algorithm is easily implemented in high level Quantitative Programming Environments (QPEs)
such as {\it Mathematica\/}, MatLab and R.
In order to obtain reasonable speed,
previous ARMA maximum likelihood algorithms are usually implemented in C or some other machine efficient language.
With our algorithm it is easy to do maximum likelihood estimation
for long time series directly in the QPE of your choice.
The new algorithm is extended to obtain the MLE for the mean parameter.
Simulation experiments which illustrate the effectiveness of the new algorithm are discussed.
{\it Mathematica\/} and R packages which implement the algorithm discussed in this
paper are available (McLeod and Zhang, 2007).
Based on these package implementations, it is expected that the interested researcher
would be able to implement this algorithm in other QPE's.


\newpage
{\noindent \bf 1. Introduction\hfill}
\bigskip

The ARMA$(p,q)$ model may be written in operator notation as
$\phi(\B)(z_t-\mu)=\theta(\B)a_t$,
where
$\B$ is the backshift operator on $t$,
$\phi(\B) = 1 -\phi_1 \B - ... - \phi_p \B^p$,
$\theta(\B) = 1 -\theta_1 \B - ... - \theta_p \B^q$, $\mu$ is the mean of $z_t$
and $a_t$ is assumed to be Gaussian white noise with mean zero and variance $\sigma_a^2$.
It is assumed that $z_t$ is causal-stationary and invertible so that all roots of
$\phi(\B)\theta(\B)=0$ are outside the unit circle.
For model identifiability it is assumed that $\phi(\B)$ and $\theta(\B)$ have no common factors.
Given $n$ consecutive observations from this time series model,
$z_1,\ldots,z_n$, the log-likelihood function was discussed by Box, Jenkins and Reinsel (1994),
as well as many other authors.
Other asymptotically first-order efficient methods are available,
such as the {\it HR\/}\ algorithm
(Hannan and Rissanen, 1982) but many researchers prefer methods of estimation and inference
based on the likelihood function
(Barnard, Jenkins and Winsten, 1962; Fisher, 1973; Box and Luce\~no, 1997, \S 12B)
and Taniguchi (1983) has shown that MLE is second-order efficient.
Some of the widely used algorithms for ARMA likelihood evaluation
are listed in Box and Luce\~no (1997, \S 12B).
All of these algorithms require $O(n)$ flops per likelihood evaluation.
The algorithm presented in \S 3 requires only $O(1)$ flops per evaluation and
so is much more efficient for longer time series.
This is especially important when implementing the algorithm in a high level QPE.
For example, one may be interested in forecasting long time series in biomedical
signal processing using MatLab (Baura, 2002, \S 7.1).
In \S 2 we discuss the AR$(p)$ case and in \S 3 the extension to the ARMA$(p,q)$ case.

\bigskip
{\noindent \bf 2.  AR$(p)$ Case\hfill}
\bigskip

\bigskip
{\noindent \it 2.1. Exact Likelihood Function\hfill}
\bigskip

It follows from Champernowne (1948, eq. 3.5) and Box, Jenkins and Reinsel (1994, eqn. A7.4.10)
that the log-likelihood function may be written
\begin{equation}
L(\phi,\mu,\sigma_a^2)=-{n\over 2}\log(\sigma_a^2)-{1\over 2}\log(g_p)-S(\phi,\mu)/(2 \sigma_a^2),
\newcounter{loglikelihood}
\setcounter{loglikelihood}{\value{equation}}
\end{equation}
where $\phi=(\phi_1,\ldots,\phi_p)$,
$g_p=\det(\Gamma_n \sigma_a^{-2})=\det(\Gamma_p \sigma_a^{-2})$, $\Gamma_n$ is the covariance matrix of
$n$ successive observations,
\begin{equation}
S(\phi,\mu)=\beta^{\prime} D \beta,
\newcounter{SumOfSquares}
\setcounter{SumOfSquares}{\value{equation}}
\end{equation}
where $D$, the Champernowne matrix, is the $(p+1) \times (p+1)$ matrix with $(i,j)$-entry,
\begin{equation}
D_{i,j}=D_{j,i}=(z_i-\mu) (z_j-\mu) + \ldots + (z_{n+1-j}-\mu) (z_{n+1-i}-\mu)
\newcounter{Dij}
\setcounter{Dij}{\value{equation}}
\end{equation}
and $\beta =(-1,\phi)$.
It should be pointed out that Champernowne (1948, p.206) assumes $n>2p$.
However, it may be shown (McLeod and Zhang, 2007) that
eqn. (\theSumOfSquares) is valid if and only if $n \ge 2p$.

Maximizing over $\sigma_a^2$, the concentrated log-likelihood may be written
\begin{equation}
L_c(\phi, \mu)=-{n\over 2}\log(S(\phi, \mu)/n)-{1\over 2}\log(g_p).
\newcounter{loglikelihoodM}
\setcounter{loglikelihoodM}{\value{equation}}
\end{equation}
As in Jones (1980),
the parametrization using partial autocorrelations
(Barndorff-Nielsen and Schou, 1973),
\begin{equation}
(\phi_1,\ldots,\phi_p) \longleftrightarrow (\zeta_1,\ldots,\zeta_p)
\newcounter{bijection}
\setcounter{bijection}{\value{equation}}
\end{equation}
may be used to constrain the optimization.
In the reparameterized model,
\begin{equation}
g_p=\prod_{j=1}^p (1-\zeta_j^2)^{-j}.
\newcounter{CovarianceDet}
\setcounter{CovarianceDet}{\value{equation}}
\end{equation}
The Burg estimators are used as initial estimates since
they are more accurate than the Yule-Walker estimates in many situations
(Percival and Walden, 1993, p.414; Zhang and McLeod, 2006b).
Like the Yule-Walker estimates, the Burg estimates are always inside the admissible region
and may be efficiently computed using the Durbin-Levinsion recursion (Percival and Walden, 1993, p.452).
Modern QPEs provide various built-in algorithms for nonlinear function optimization
which may be used to obtain the MLE of $\phi$.
Since the sample mean, $\bar z = (z_1 + \cdots +z_n)/n$,
is an asymptotically fully efficient estimate of $\mu$, it is often used in place of the MLE.
This algorithm using the sample mean to estimate $\mu$ and then MLE for the other
parameters will be denoted by {\it SampleMean} in the following sections.

If the sample mean is used, $\mu$ may be replaced by $\bar z$ in (\theDij) and so
after the initial evaluation, repeated evaluations of (\theloglikelihoodM)
require $O(1)$ flops, which explains why the new algorithm is efficient for long time series.
Since it practice $p$ is considered fixed, it is not included in the asymptotic flop count.

\bigskip
{\noindent \it 2.2. Exact MLE for the Mean Parameter\hfill}
\bigskip

The exact MLE for the mean may be obtained by simply optimizing the log-likelihood
function given in (\theloglikelihoodM).
However, this would then require $O(n)$ flops per function evaluation.
A more efficient approach is now presented.

Assuming that $\phi$ is known,
the exact MLE is given by,
\begin{equation}
\hat \mu = {1_n^\prime \Gamma_n^{-1} z \over {1_n^\prime \Gamma_n^{-1} 1_n}},
\newcounter{ExactMu}
\setcounter{ExactMu}{\value{equation}}
\end{equation}
where $1_n$ denotes the $n$ dimensional column vector with all entries equal to $1$,
$1_n^\prime$ denotes its transpose
and $z=(z_1,\ldots,z_n)$.
Since $\hat \mu$ does not depend on $\sigma_a^2$, we may assume without loss of generality
that $\sigma_a^2=1$.
Direct evaluation of (\theExactMu) using the exact inverse matrix derived by
Siddiqui (1958) would require $O(n^2)$ flops.
A more efficient approach may be developed using the inverse matrix
result of Zinde-Walsh (1988).
Zinde-Walsh (1988, eqn. 3.2) showed that
\begin{equation}
\Gamma_n^{-1}=\dot \Gamma_n-\Omega,
\newcounter{ZindeWalsh}
\setcounter{ZindeWalsh}{\value{equation}}
\end{equation}
where $\dot \Gamma_n$ denotes the $n\times n$ matrix with $(i,j)$-entry given by
$\gamma_{i-j}^{(u)}$, where
$\gamma_k^{(u)}=\Cov(u_t, u_{t-k})$, $u_t = \phi(\B) a_t$
and $\Omega$ is a zero matrix except for $p\times p$ submatrices in the upper-left and lower-right corners.
The $(i,j)$-entry of the submatrix of $\Omega$ in the upper-left corner is
\begin{equation}
\Omega_{i,j}=\sum_{k=\min(i,j)}^{p-|i-j|} \phi_k \phi_{k+|i-j|}.
\newcounter{Omegaij}
\setcounter{Omegaij}{\value{equation}}
\end{equation}
The matrix in the lower-right corner is just the transpose of the upper-left corner submatrix.
Using the above results it was found that,
\begin{eqnarray}
1_n^\prime \Gamma_n^{-1} &=& 1_n^\prime \phi^2(1)
-(\epsilon_1,\ldots,\epsilon_p,0,\ldots,0,\epsilon_p,\ldots,\epsilon_1)
\nonumber\\
&&-(\kappa_1,\ldots,\kappa_p,0,\ldots,0,\kappa_p,\ldots,\kappa_1),
\newcounter{NewResult}
\setcounter{NewResult}{\value{equation}}
\end{eqnarray}
where
$\phi(1)=1-\phi_1-\ldots-\phi_p$,
$\epsilon=1_n^{\prime} \Omega$,
\begin{equation}
\epsilon=(\epsilon_1,\ldots,\epsilon_p,0,\ldots,0,\epsilon_p,\ldots,\epsilon_1)
\newcounter{epsilons}
\setcounter{epsilons}{\value{equation}}
\end{equation}
and
\begin{equation}
\kappa_i=\sum_{k=1}^{i} \gamma_k^{(u)}.
\newcounter{kappaij}
\setcounter{kappaij}{\value{equation}}
\end{equation}
Using (\theNewResult), $\hat \mu$ can now be evaluated in $O(n)$ flops.
Note that this evaluation will only typically be two or three times in the full MLE
algorithm outlined below.

An iterative algorithm, {\it MeanMLE}, is used for the simultaneous joint MLE
of $(\phi_1,\ldots,\phi_p,\mu)$,

\begin{description}
\item[Step 0] Set the maximum number of iterations, $M \leftarrow 5$.
Set the iteration counter, $i\leftarrow 0$.
Set $\hat \mu^{(0)}\leftarrow \bar z$, where $\bar z$ is the sample mean.
Obtain initial parameter values $\hat \phi^{(0)}_k,\ k=1,\ldots,p$ using the Burg
algorithm or set $\hat \phi^{(0)}_k = 0,\ k=1,\ldots,p$.
Set $\ell_0 = L_c(\hat \phi^{(0)}, \hat \mu^{(0)})$.
\item[Step 1] Obtain $\hat \phi^{(i+1)}_k,\ k=1,\ldots,p$ by numerically
maximizing $L_c(\phi,\hat \mu^{(i)})$ over $\phi$.
Set $\ell_{i+1} = L_c(\hat \phi^{(i+1)}, \hat \mu^{(i)})$.
\item[Step 2] Using $\hat \phi^{(i+1)}$ evaluate $\hat \mu^{(i+1)}$.
\item[Step 3] Terminate when $\ell_{i+1}$ has converged or $i>M$.
Otherwise set $i\leftarrow i+1$ and
return to Step 1 to perform the next iteration.
\end{description}

Convergence usually occurs in two or three iterations.

\bigskip
{\noindent \it 2.3. Champernowne Matrix Computation\hfill}
\bigskip

$D_{i,j}$ has $n-(i+1)-(j+1)$ terms so each term requires $O(n)$ flops.
If the sample mean is used, this computation only has to be done once,
but if the exact MLE
for the mean is used, $D$ must be computed several times.
It may be shown that $D=C-E$, where the $(i,j)$-entry of the matrix $C$ may be written,
$C_{|i-j|}$, where $C_k = z_1 z_k + \ldots + z_{n-k} z_n$.
The $(i,j)$-entry for the matrix $E$ may be computed sequentially
$E_{i+1,j+1} \linebreak[0]= \linebreak[0]E_{i,j}+z_i z_j +z_{n+1-i} z_{n+1-j}$, $i<j$.
Using the above results reduces the flop count for the matrix $D$ slightly.

\bigskip
\strut
{\noindent \bf 3. ARMA Maximum Likelihood Estimation\hfill}
\bigskip

Previous AR-approximation methods for fitting MA$(q)$ and ARMA$(p,q)$ were based on first fitting a
suitable high-order autoregressive approximation
(Durbin, 1959; Parzen, 1969; Hannan and Rissanen, 1982; Wahlberg, 1989; Choi, 1992 \S 4.1).
The next step is to use the fitted $\AR$ model to estimate an MA$(q)$ or ARMA$(p,q)$ model.
As noted by McClave (1973), this approach can lead to biased estimates which have
larger mean-square error than the MLE.

Instead of directly fitting an autoregressive model to the time series,
our new method is based on approximating the exact likelihood function for
the ARMA$(p,q)$ model by the likelihood function for a suitable
high-order autoregression.
The approximating autoregression of order $r$ is determined as the
minimum mean-square error (MMSE) linear predictor of order $r$ for the ARMA$(p,q)$ model,
$\varphi(B)(z_t - \mu) = a_t$, where
$\varphi(\B)=1-\varphi_1 \B - \ldots - \varphi_r \B^r$.
By taking $r$ sufficiently large, an accurate approximation to the exact ARMA$(p,q)$
likelihood may be obtained.
In practice $r=30$ is sufficient for many ARMA models as we will now show.

The Kullback-Leibler discrepancy may be used to choose a suitable $r$.
Letting $\Sigma_{\phi,\theta}$ and $\Sigma_{\varphi}$ denote the covariance
matrices for the ARMA$(p,q)$ and its AR$(r)$ approximation,
the Kullback-Leibler discrepancy may be written (Ullah, 2002, eqn. 5),
\begin{equation}
{\cal I} = {1\over2}( \tr \Sigma_{\phi,\theta} \Sigma_{\varphi}^{-1} - \log |\Sigma_{\phi,\theta}|/|\Sigma_{\varphi}| -n).
\newcounter{KLDiscrepancy}
\setcounter{KLDiscrepancy}{\value{equation}}
\end{equation}
Figure 1 displays a plot of $\cal I$ in the case of an $\MA(1)$ model with
$\theta_1=0.9$ and $n=200$.
It is seen that $r=30$ works well even for this model with a parameter near
the non-invertible boundary.
It appears that $r=30$ is adequate for many sorts of models occurring in applications although
as the parameters move very close to the non-invertible boundary,
our approximates requires larger $r$ and fails entirely when the boundary is reached.
A {\it Mathematica\/} notebook to compute and plot the Kullback-Leibler discrepancy
for the ARMA$(p,q)$ and its AR$(r)$ approximation is available (McLeod and Zhang, 2007).

\medskip
\centerline{[Figure 1 here]}
\medskip

In practice, as shown by simulation in \S 4.2,
our method with $r=30$ can still be used even when there is a root on the boundary
but the statistical efficiency relative to existing exact MLE algorithms is reduced.
Models with a root on the non-invertible boundary usually indicate over-differencing and
may be avoided by refitting with an alternative model specification (Zhang and McLeod, 2006a).

After a suitable $r$ has been chosen, the ARMA likelihood may be
obtained from (\theloglikelihoodM),
\begin{equation}
L_c(\phi, \theta, \mu) = L_c(\varphi, \mu),
\newcounter{loglikelihoodARMA}
\setcounter{loglikelihoodARMA}{\value{equation}}
\end{equation}
where $\varphi = (\varphi_1, \ldots, \varphi_r)$.
Then $L_c(\phi, \theta, \mu)$ may be maximized using a built-in optimization function.
The algorithm given in \S 2.2 may be used to compute the exact MLE for the mean by
using this AR$(r)$ approximation.
As shown in \S 4.3, this algorithm works as well as existing exact MLE algorithms
for the mean in ARMA$(1,1)$ models.

In {\it Mathematica\/}, MatLab and in R,
nonlinear optimization functions which can handle box constraints are available.
In this case it is useful to reparametrize the ARMA model as suggested by
Monahan (1984) using the transformation of Barndorff-Nielsen and Schou (1973).
Alternatively, if only an unconstrained optimization function is available then
a penalty function approach may be used to constrain the parameters to the
admissible region.
This penalty function approach has been used for many years
with the Powell (1964) algorithm in our MHTS Time Series Package
(McLeod and Hipel, 2007)
for a wide variety of MLE problems in time series analysis
(Hipel and McLeod, 1994).

Usually it is most expedient to set the initial parameter estimates to zero.
In case of difficulty with convergence, initial estimates may be obtained
(Hannan and Rissanen, 1982)
by fitting a high order autoregression to provide estimates of the innovations
and then using linear regression to estimate the parameters $\phi$ and $\theta$.
Experience suggests, as is illustrated in \S 4.1, computing initial
parameter estimates in the ARMA case usually does not significantly increase the speed and,
in practice, convergence is rarely an issue.
In particular, convergence was obtained for all models fitted in \S 4 without difficulty.

A simple alternative to the MMSE linear predictor approximation is to just use
the truncated inverted form of model (Box, Jenkins and Reinsel, 1994, \S 4.2.3),
$\pi(B) (z_t-\mu)=a_t$, where
$\pi(B)=1-\pi_1 B - \ldots - \pi_r B^r$.
The coefficients $\pi_k$, $k=1,\ldots,r$
are obtained from $\pi_k =\phi_k + \theta_1 \pi_{k-1} - \ldots - \theta_q \pi_q - \phi_k$
using boundary conditions
$\pi_0 = 1; \pi_k = 0\ {\rm if}\ k<0$
and $\phi_k = 0\ {\rm if}\ k>p.$
When $r$ is chosen large enough, this approximates the MMSE predictor (Brockwell and Davis, 1991, \S 5).
However, for fixed $r$ there will always be parameter values in the admissible
ARMA$(p,q)$ region for which $\varphi(B)=0$ has roots outside the admissible region
for a causal-stationary AR$(r)$.
As shown in Table 1, the MMSE predictor provides a much more accurate approximation
in terms of the Kullback-Leibler discrepancy.
For these reasons the MMSE linear predictor approximation is used.

\medskip
\centerline{[Table 1  here]}
\medskip

\bigskip
{\noindent \bf 4. Illustrative Examples\hfill}
\medskip

The primary purpose of the illustrative examples presented in this section is to
demonstrate the usefulness of our algorithm and correctness
of our implementations in {\tty R} and {\it Mathematica\/}.
For this purpose, our algorithm is also compared with existing MLE
algorithms.

\bigskip
{\noindent \it 4.1. Timings \hfill}
\bigskip

Timings for the algorithms described in \S 3
were obtained in {\it Mathematica\/} and R on a Windows XP PC Pentium 4.
The ARMA$(1,1)$ model with $\phi_1=0.9$ and $\theta=0.5$ was selected as
typical of order $(1,1)$ models which might occur in practice.
This model was simulated 25 times for series of length $n=10^k,\ k=2,3,\ldots,6$
and the average time needed for fitting the model was determined.
Timings were also compared to {\it HR\/} (Hannan and Rissanen, 1982).
The {\it HR\/} algorithm does not require non-linear optimization and only requires
linear least squares and residual computation.
The built-in least squares algorithms in {\it Mathematica} and {\tty R} were
used.
The effects of initial values and MLE estimation of the mean were also examined.
The initial value options also examined were {\it Origin}, {\it XInit} and {\it HRInit}
corresponding respectively to initializing the nonlinear optimization algorithm at $0.0$ for all parameter values
except the mean, using exact known parameter values or using the Hannan-Rissanen
estimates as initial parameter settings.
The algorithms for estimating the mean, {\it SampleMean} and {\it MeanMLE\/}, are also compared.
The {\it MeanMLE} refers to the algorithm in \S 2.2 and
{\it SampleMean} to just using $\bar z$ as in \S 2.1.
In the R timings we also compared our algorithms with the
built-in {\tty R} algorithms {\tty arima} and {\tty arima0}.
These algorithms implement the state-space Kalman filter algorithm
given in Durbin and Koopman (2001).
Further details of this implementation (Ripley, 2002) indicate that
this algorithm is coded in C and then interfaced to {\tty R}.

\medskip
\centerline{[Table 2 about here]}
\medskip

Comparing {\it Origin\/} with {\it HR\/}, our algorithm is much faster for larger $n$.
Although {\it HR\/} is faster than {\it Origin\/} for small $n$ this is probably not important
since both algorithms are very fast and {\it Origin\/} which uses the MLE method
is preferred anyway -- especially for small $n$.
Since the computing time required by {\it HRInit\/} does not include the initialization times
needed by {\it HR\/} itself, it is clear from Table 2, that if these are added to
{\it HRInit\/}, the initialization is normally not worthwhile in terms of reducing
computer time.
Even with {\it XInit} when the exact initial values are used, this only results in a modest improvement in speed.
It is seen that in terms of speed {\it Mathematica\/} outperforms {\tty R} except when $n$ is very large.
These timings also demonstrate that the {\it Mathematica} and {\tty R} implementations
of our algorithms are suitable for even very large $n$.
Given the high-overhead imposed by the interpretive {\tty R} language,
the performance of our algorithms is not unreasonable in practice even though
in most cases it is slower than {\tty arima} and {\tty arima0}.

\bigskip
{\noindent \it 4.2. Comparison with Durbin's Algorithm \hfill}
\bigskip

The statistical efficiency of {\it Durbin\/}, the algorithm of Durbin (1959) for MA$(q)$
estimation, is compared with {\it SampleMean\/} and exact MLE as implemented in {\tty R} in
{\tty arima}.
For each parameter value $\theta_1 = 0, \pm 0.3, \pm 0.5, \pm 0.9, \pm 1$,
and for each series length $n=50,100,200,400$
one thousand time series were simulated.
The empirical statistical efficiency may be taken as the empirical MSE of the exact MLE algorithm divided
by the empirical MSE of {\it SampleMean\/}.
Similarly, for the efficiency for the {\it Durbin\/} algorithm.
The variance of the estimated efficiency may be derived using a Taylor
series linearization.
Details of this derivation as well as a comparison with the bootstrap
variance estimate are given in our online supplement (McLeod and Zhang, 2007).
In Figure 2, a trellis plot compares these efficiencies.
In each plot, the vertical line running through the plotted point indicates a
95\% confidence interval for that efficiency.
From this plot, we see that {\it SampleMean\/} has efficiency very close to 1 except
when the parameter $\theta_1 = \pm 1$ when it is less efficient and when $\theta_1 = 0.9$
it is super-efficient.
In the super-efficiency cases, the efficiency approaches 1 as $n$ increases.
The efficiency of {\it Durbin\/} is generally much less than {\it SampleMean\/} but it approaches 1
as $n$ gets larger provided the parameter is not on the boundary.

The results shown in Figure 2 were replicated using our {\it Mathematica} implementation
of {\it SampleMean} and
the exact MLE algorithm for the $\MA(1)$ given in McLeod and Quenneville (2001).

\medskip
\centerline{[Figure 2 here]}
\medskip

\bigskip
{\noindent \it 4.3. Finite Sample Efficiency of the Sample Mean \hfill}
\bigskip

If the parameters $\phi_1,...,\phi_p,\theta_1,\ldots,\theta_q$ are known,
the exact MLE for the mean is given by eqn. (\theExactMu).
It is also the best linear unbiased estimate {\it BLUE}.
Another estimate of $\mu$ is simply the sample mean,
$\bar z = (z_1+\ldots+z_n)/n$.
The exact efficiency for $\bar z$ vs. the {\it BLUE} for a series of length $n$ may be written,
\begin{equation}
{\cal E} = n^2 / ( (1_n^{\prime} \Gamma 1_n) (1_n^{\prime} \Gamma^{-1} 1_n) ).
\newcounter{EffSampleMean}
\setcounter{EffSampleMean}{\value{equation}}
\end{equation}
In actual applications, the ARMA parameters are not known.
In our simulation study, we compare two MLE methods
for estimating the mean.
The MLE methods are the {\it MeanMLE} algorithm of \S 2.2 and the {\tty R} function
{\tty arima}.
With each of these MLE methods, the empirical efficiency of $\bar z$ vs. the MLE
estimate of $\mu$ based on $10^3$ simulations for series of lengths $n=50,100,200$
for the ARMA$(1,1)$ model at each parameter setting.
These empirical efficiencies are compared with the exact efficiency of $\bar z$ vs.
{\it BLUE} given in eqn. (\theEffSampleMean) and all results are displayed in Table 3.
Both {\it MeanMLE} and {\tty arima} are closely efficient and there is general agreement
with the {\it BLUE} except when $\phi_1 = 0$ and $\theta_1=0.9, 0.95$.
The simulation experiment confirms that {\tty MeanMLE} is working correctly
as expected and this was its main purpose.

Since the sample mean is asymptotically efficient in ARMA$(p,q)$ models (Brockwell and Davis, 1991, \S 7.1)
it would be expected the efficiencies would get closer to 1 as $n$ increases and
it is seen that in many cases this holds.
However it is surprising that even for $n=200$, some sample efficiencies
are quite low for both {\it MeanMLE} and {\tty arima}.
This fact does not previously appear to have been observed in the ARMA case
although Samarov and Taqqu (1988) found asymptotic inefficiency in a situation
which we will now discuss briefly.

It should be noted that the ARMA models where the sample mean efficiency
is low have an extremely high frequency spectrum.
The spectral density and autocorrelation plots of the models in Table 3 are given
in McLeod and Zhang (2007).
The models for which the sample mean is inefficient all have strong
negative autocorrelation but are better characterized in terms of
the spectral density function.
All models for which the sample mean efficiency is less than 10\% efficient are all
characterized by a high frequency spectrum in which the high frequencies are
more than one hundred times the power of the low frequencies, that is,
the ratio of the spectral density evaluated at the Nyquist frequency divided
by the spectral density evaluated at the origin is larger than 100.
This situation may be called, infrared-catastrophe since it
seems unrealistic in any time series applications with actual scientific data.

Previously Samarov and Taqqu (1988) showed that asymptotically the sample mean
can be very inefficient for hyperbolic decay time series (McLeod, 1988)
in the antipersistent case which corresponds to the infrared-catastrophe case for these models.
In all other hyperbolic-decay cases,
including the fractional ARMA case in eqn. (16),
the asymptotic efficiency is above 98\% (Samarov and Taqqu, 1988, Table 1).

This simulation experiment was repeated using the {\it SampleMean\/} algorithm implemented
in {\it Mathematica\/} and similar results were obtained (McLeod and Zhang, 2007).

\medskip
\centerline{[Table 3  here]}
\medskip

\bigskip
{\noindent \bf 5. Conclusion\hfill}
\medskip

{\it Mathematica\/} and {\tty R}\ packages that implement the ARMA maximum likelihood
algorithms described in this paper are available (McLeod and Zhang, 2007).
In addition simulation scripts to obtain the results reported in this article
are also available so the interested can easily reproduce and/or extend our simulation
results using either {\it Mathematica\/} or {\tty R}.

Our algorithms are suitable for use with long time series.
But the principal advantage of our algorithms for maximum likelihood estimation of ARMA
models is that they may easily be implemented directly in high-level QPEs.
Using the {\tty R}\ and {\it Mathematica} packages, it is relatively straightforward to implement ARMA maximum likelihood in other high level QPEs.
QPEs such as MatLab and Strata as well as {\tty R} and {\it Mathematica} are becoming
important in teaching statistical methods so it is expected our algorithm
will be useful teaching time series analysis in such computing environments.

The AR-likelihood approximation technique of this paper could be used for other types
of linear time series models.
It would be relatively straightforward to extend the methods of this paper to
multiplicative seasonal and subset ARMA models.
It may also be possible to develop an extension to the vector ARMA models case.
Another interesting family of linear time series models
are the fractional ARMA time series (Hipel and McLeod, 1994, Ch. 11;
Brockwell and Davis, \S 13.2) defined by
\begin{equation}
\phi(\B) \nabla^d (z_t-\mu)=\theta(\B)a_t,
\newcounter{FractionalARMA}
\setcounter{FractionalARMA}{\value{equation}}
\end{equation}
where $d \in (-0.5,0.5)$.
Figure 3 shows the Kullback-Leibler discrepancy, ${\cal I}$ for the case
of fractionally differenced white noise, $p=0$ and $q=0$, with long-memory
parameter $d=0.1,0.2,0.3,0.4$.
When $d \in (0,0.2)$, $r=30$ is adequate but much higher orders
may be needed for more strongly persistent time series such as when $d \ge 0.4$.
In the case such strongly persistent time series our suggested AR approximation
may not be useful.

\medskip
\centerline{[Figure 3 here]}
\medskip

{\noindent \bf Acknowledgements\hfill}
Both authors were supported by NSERC Discovery Grants.
The authors would like to thank the referees for helpful comments and
their careful reading of our work.

\bigskip

\medskip
\newpage
\noindent
{\bf Referemces\hfill}
\parindent 0pt
\medskip

\hind

\hind
Barnard, G.A., Jenkins, G.M. and Winston, C.B. 1962.
Likelihood inference and time series.
{\it Journal of the Royal Statistical Society, B} 125, 321--372.

\hind
Barndorff-Nielsen, O. and Schou, G., 1973.
On the parametrization of autoregressive models by partial autocorrelations.
{\it Journal of Multivariate Analysis} 3, 408--419.

\hind
Baura, G.D., 2002.
{\it System Theory and Practical Applications of Biomedical Signals.}
Wiley, New York.


\hind
Box, G.E.P., Jenkins, G.M. and Reinsel, G.C.,  1994.
{\it Time Series Analysis: Forecasting and Control\/}.
3rd Ed., Holden-Day, San Francisco.

\hind
Box, G.E.P. and Luce\~no, A., 1997.
{\it Statistical Control by Monitoring and Feedback Adjustment\/}.
Wiley, New York.

\hind
Brockwell, P.J. and Davis, R.A., 1991.
{\it Time Series: Theory and Methods\/}.
(2nd edn.) Springer-Verlag,  New York.

\hind
Champernowne, D.G., 1948.
Sampling theory applied to autoregressive sequences.
{\it Journal of the Royal Statistical Society B} 10, 204--242.

\hind
Choi, B., 1992.
{\it ARMA Model Identification\/}.
Springer-Verlag, New York.

\hind
Durbin, J., 1959.
Efficient estimation of parameters in moving-average models.
{\it Biometrika\/} 46, 306--316.

\hind
Durbin, J. and Koopman, S.J., 2001.
{\it Time Series Analysis by State Space Methods.}
Oxford University Press, Oxford.

\hind
Fisher, R.A. 1973.
{\it Statistical methods and scientific inference.}
Hafner Press, New York.

\hind
Hannan, E.J. and Rissanen, J., 1982.
Recursive estimation of mixed autoregressive-moving average order.
{\it Biometrika\/} 69, 81--94.

\hind
Hipel, K.W. and McLeod, A.I., 1994.
{\it Time Series Modelling of Water Resources and Environmental Systems.}
Elsevier, Amesterdam.
Reprint, www {\tty http://www.stats.uwo.ca\linebreak[0]/faculty/\linebreak[0]aim/\linebreak[0]1994Book/}.

\hind
Jones, R.H., 1980.
Maximum likelihood fitting of ARMA models to time series with missing observations.
{\it Technometrics\/} 22, 389--395.

\hind
McClave, E.J., 1973.
On the bias of AR approximation to moving averages.
{\it Biometrika\/} 60, 599-605.

\hind
McLeod, A.I., 1998.
Hyperbolic decay time series.
{\it Journal of Time Series Analysis\/} 19, 473--484.

\hind
McLeod, A.I. and Quenneville, B., 2001.
Mean likelihood estimators.
{\it Statistics and Computing\/} 11, 57--65.

\hind
McLeod, A.I. and Hipel, K.W., 2007.
McLeod-Hipel Time Series Package,
(www {\tty http://www.stats.uwo.ca/faculty/aim/epubs/mhts/}).

\hind
McLeod, A.I. and Zhang, Y., 2007.
Online supplements to ``Faster ARMA Maximum Likelihood Estimation'',
www ({\tty http://www.stats.uwo.ca/faculty/aim/2007/faster/}).

\hind
Monahan, J.F., 1984.
A note on enforcing stationarity in autoregressive-moving average models.
{\it Biometrika} 71, 403--404.

\hind
Parzen, E., 1969.
Multiple time series modeling.
In {\it Multivariate Analysis II\/} ed. P. Krishnaiah, 389-409.
Academic Press, New York.

\hind
Percival, D.B. and Walden, A.T., 1993.
{\it Spectral Analysis for Physical Applications\/}.
Cambridge University Press, Cambridge.

\hind
Powell, M.J.D., 1964.
An efficient method for finding the minimum of a function of several
variables without calculating derivatives.
{\it Computer Journal\/} 7, 155--162.

\hind
Ripley, B.D., 2002.
Time Series in R.
{\it R News} 2, 2--7.

\hind
Samarov, A. and Taqqu, M., 1988.
On the efficiency of the sample mean in long memory noise.
{\it Journal of Time Series Analysis\/} 9, 191--200.

\hind
Siddiqui, M.M., 1958.
On the inversion of the sample covariance matrix in a stationary autoregressive process.
{\it Annals of Mathematical Statistics} 29, 585--588.




\hind
Taniguchi, M., 1983.
On the second order asymptotic efficiency of estimators of gaussian ARMA processes.
{\it The Annals of Statistics} 11, 157--169.

\hind
Ullah, A., 2002.
Use of entropy and divergence measures for evaluating econometric approximations
and inference.
{\it Journal of Econometrics\/} 107, 313--326.

\hind
Wahlberg, B., 1989.
Estimation of autoregressive moving-average models via high-order autoregressive
approximations.
{\it Journal of Time Series Analysis\/} 10, 283--299.

Zhang, Y. and McLeod, A.I.,  2006a.
Fitting MA$(q)$ models in the closed invertible region.
{\it Statistics and Probablity Letters\/} 76, 1331--1334.

\hind
Zhang, Y. and McLeod, A.I., 2006b.
Computer algebra derivation of the bias of Burg estimators.
{\it Journal of Time Series Analysis\/} 27, 157--165.

\hind
Zinde-Walsh, V., 1988.
Some exact formulae for autoregressive moving average processes.
{\it Econometric Theory} 4, 384--402.

\clearpage
\newpage
\begin{table}
\caption{
Kullback-Leibler discrepancy for AR$(r)$ approximation to a $\MA(1)$
with $\theta_1 = 0.95$ and $n=200$ using the MMSE approximation and
the approximation based on truncating the inverted form of the model.
}

\begin{tabular}{lll}
\cr
\hline
\noalign{\smallskip}
$n$&MMSE&Truncated \cr
\cline{2-3}
\noalign{\smallskip}
10  &4.99   &31.20\cr
20  &1.23   &10.70\cr
30  &0.38   &3.77\cr
40  &0.12   &1.42\cr
50  &0.04   &0.62\cr
\noalign{\smallskip}
\noalign{\hrule}
\end{tabular}
\end{table}

\clearpage
\newpage
\begin{table}
\caption{
Average CPU time in seconds with {\tty R} and {\it Mathematica} for fitting
the ARMA$(1,1)$ model with $\phi_1=0.9$ and $\theta=0.5$ using {\it SampleMean\/},
{\it MeanMLE\/} and the Hannan-Rissanen estimator.
In the {\tty R}\ case, built-in functions
{\tty arima} and {\tty arima0} are also used.
Twenty-five replications for series of length $n=10^k,\ k=2,3,\ldots,6$ were done.
The case where the mean is estimated by the sample average is compared
with the MLE for each algorithm.
The effect of initial parameter settings is also examined.
The settings {\it Origin}, {\it XInit} and {\it HRInit} correspond
to setting $(\phi_1,\theta_1)$ equal to $(0,0)$, $(0.9,0.5)$ or using the estimator of
Hannan-Rissanen respectively.
}
\begin{tabular}{llllll}
\cr
\hline
\noalign{\smallskip}
&\multicolumn{5}{c}{$n$} \cr
method&$10^2$&$10^3$&$10^4$&$10^5$&$10^6$ \cr
\cline{2-6}
\noalign{\smallskip}
&\multicolumn{5}{c}{Timings in {\tty R}} \hfil \cr
\noalign{\smallskip}
\cline{2-6}
\noalign{\smallskip}
&\multicolumn{5}{c}{\it SampleMean} \cr
\noalign{\smallskip}
Origin    & 0.47&0.47&0.70&2.13&8.63\cr
HR        &0.24&1.05&10.2&92.0&902.\cr
XInit     &0.32&0.27&0.46&1.40&6.81\cr
HRInit    &0.29&0.27&0.57&1.70&7.78\cr
arima     &0.02&0.04&0.24&1.85&13.6\cr
arima0    &0.01&0.01&0.02&0.18&1.78\cr
&\multicolumn{5}{c}{\it MeanMLE} \cr
Origin    &1.00&0.85&1.03&3.14&16.70\cr
XInit     &0.90&0.68&0.82&2.80&17.15\cr
HRInit    &0.80&0.58&0.89&2.91&16.59\cr
arima     &0.05&0.19&0.74&3.94&32.54\cr
arima0    &0.03&0.03&0.07&0.88&13.62\cr
\noalign{\smallskip}
\cline{2-6}
\noalign{\smallskip}
&\multicolumn{5}{c}{Timings in {\it Mathematica}} \cr
\noalign{\smallskip}
\cline{2-6}
\noalign{\smallskip}
&\multicolumn{5}{c}{\it SampleMean} \cr
\noalign{\smallskip}
\cline{2-6}
Origin      &0.26     &0.29     &0.36    &0.90    &4.97\cr
HR          &0.01     &0.05     &0.54    &5.00    &47.4\cr
XInit       &0.27     &0.30     &0.36    &0.89    &4.96\cr
HRInit      &0.27     &0.29     &0.35    &0.88    &4.99\cr
&\multicolumn{5}{c}{\it MeanMLE} \cr
Origin      &0.75     &0.89     &1.14    &4.10    &30.46\cr
XInit       &0.78     &0.90     &1.14    &4.11    &30.71\cr
HRInit      &0.75     &0.89     &1.10    &4.10    &30.49\cr
\noalign{\smallskip}
\cline{2-6}
\end{tabular}
\end{table}

\vfill\eject
\clearpage
\newpage
\begin{table}
\caption{
Empirical efficiency of the sample mean vs.
three other methods: BLUE, {\it MeanMLE} and {\tty arima}.
Each empirical efficiency is based on 1000 simulations
for ARMA$(1,1)$ models with $n=50,100,200$.
}

\begin{tabular}{llllllllll}
\cr
\hline
\noalign{\smallskip}
& & \multicolumn{7}{c}{$\theta_1$} \cr
\cline{4-10}
\noalign{\smallskip}
$\phi_1$&algorithm&$n$&$-0.95$&$-0.9$&$-0.5$&$0.$&$0.5$&$0.9$&$0.95$\cr
\noalign{\smallskip}
\noalign{\hrule}
\noalign{\smallskip}
$-0.95$&{\it BLUE}&$50$&$1.00$&$1.00$&$0.97$&$0.75$&$0.25$&$0.01$&$0.01$
\cr
$-0.95$&{\it MeanMLE}&$50$&$1.00$&$1.00$&$0.99$&$0.74$&$0.39$&$0.02$&$0.01$
\cr
$-0.95$&{\tty arima}&$50$&$1.02$&$1.02$&$1.00$&$0.74$&$0.40$&$0.02$&$0.01$
\cr\noalign{\smallskip}
$-0.95$&{\it BLUE}&$100$&$1.00$&$1.00$&$0.98$&$0.85$&$0.38$&$0.02$&$0.01$
\cr
$-0.95$&{\it MeanMLE}&$100$&$1.00$&$1.00$&$1.00$&$0.87$&$0.58$&$0.04$&$0.01$
\cr
$-0.95$&{\tty arima}&$100$&$1.01$&$1.01$&$1.00$&$0.87$&$0.58$&$0.04$&$0.01$
\cr\noalign{\smallskip}
$-0.95$&{\it BLUE}&$200$&$1.00$&$1.00$&$0.99$&$0.92$&$0.54$&$0.03$&$0.01$
\cr
$-0.95$&{\it MeanMLE}&$200$&$1.00$&$1.00$&$0.99$&$0.91$&$0.67$&$0.06$&$0.02$
\cr
$-0.95$&{\tty arima}&$200$&$1.01$&$1.00$&$0.99$&$0.91$&$0.67$&$0.06$&$0.02$
\cr\noalign{\smallskip}
$-0.9$&{\it BLUE}&$50$&$1.00$&$1.00$&$0.99$&$0.86$&$0.39$&$0.02$&$0.01$
\cr
$-0.9$&{\it MeanMLE}&$50$&$1.00$&$1.01$&$1.00$&$0.85$&$0.59$&$0.05$&$0.02$
\cr
$-0.9$&{\tty arima}&$50$&$1.02$&$1.01$&$1.01$&$0.86$&$0.59$&$0.05$&$0.02$
\cr\noalign{\smallskip}
$-0.9$&{\it BLUE}&$100$&$1.00$&$1.00$&$0.99$&$0.92$&$0.56$&$0.03$&$0.01$
\cr
$-0.9$&{\it MeanMLE}&$100$&$1.00$&$1.00$&$1.00$&$0.95$&$0.74$&$0.07$&$0.02$
\cr
$-0.9$&{\tty arima}&$100$&$1.01$&$1.01$&$1.00$&$0.95$&$0.74$&$0.07$&$0.02$
\cr\noalign{\smallskip}
$-0.9$&{\it BLUE}&$200$&$1.00$&$1.00$&$1.00$&$0.96$&$0.71$&$0.06$&$0.02$
\cr
$-0.9$&{\it MeanMLE}&$200$&$1.00$&$1.00$&$1.00$&$0.95$&$0.81$&$0.11$&$0.03$
\cr
$-0.9$&{\tty arima}&$200$&$1.01$&$1.01$&$1.00$&$0.95$&$0.81$&$0.12$&$0.03$
\cr\noalign{\smallskip}
$-0.5$&{\it BLUE}&$50$&$1.00$&$1.00$&$1.00$&$0.99$&$0.83$&$0.13$&$0.06$
\cr
$-0.5$&{\it MeanMLE}&$50$&$1.00$&$1.00$&$1.01$&$1.00$&$0.96$&$0.26$&$0.13$
\cr
$-0.5$&{\tty arima}&$50$&$1.00$&$1.00$&$1.03$&$1.00$&$0.96$&$0.26$&$0.13$
\cr\noalign{\smallskip}
$-0.5$&{\it BLUE}&$100$&$1.00$&$1.00$&$1.00$&$0.99$&$0.91$&$0.19$&$0.07$
\cr
$-0.5$&{\it MeanMLE}&$100$&$1.00$&$1.00$&$1.00$&$1.00$&$0.95$&$0.31$&$0.13$
\cr
$-0.5$&{\tty arima}&$100$&$1.00$&$1.00$&$1.01$&$1.00$&$0.95$&$0.32$&$0.13$
\cr\noalign{\smallskip}
$-0.5$&{\it BLUE}&$200$&$1.00$&$1.00$&$1.00$&$1.00$&$0.95$&$0.30$&$0.10$
\cr
$-0.5$&{\it MeanMLE}&$200$&$1.00$&$1.00$&$1.00$&$1.00$&$0.96$&$0.43$&$0.17$
\cr
$-0.5$&{\tty arima}&$200$&$1.00$&$1.00$&$1.00$&$1.00$&$0.96$&$0.44$&$0.17$
\cr\noalign{\smallskip}
\noalign{\smallskip}
\end{tabular}
\end{table}

\vfill\eject
\clearpage
\thispagestyle{empty}
\begin{table}
\begin{tabular}{llllllllll}
\cr
\hline
\noalign{\smallskip}
& & \multicolumn{7}{c}{$\theta_1$} \cr
\cline{4-10}
\noalign{\smallskip}
$\phi_1$&algorithm&$n$&$-0.95$&$-0.9$&$-0.5$&$0.$&$0.5$&$0.9$&$0.95$\cr
\noalign{\smallskip}
\noalign{\hrule}
\noalign{\smallskip}
$0.$&{\it BLUE}&$50$&$0.99$&$0.99$&$1.00$&$1.00$&$0.96$&$0.34$&$0.18$
\cr
$0.$&{\it MeanMLE}&$50$&$1.02$&$1.02$&$1.02$&$1.03$&$1.03$&$1.03$&$1.03$
\cr
$0.$&{\tty arima}&$50$&$1.03$&$1.03$&$1.03$&$1.02$&$1.03$&$1.03$&$1.03$
\cr\noalign{\smallskip}
$0.$&{\it BLUE}&$100$&$1.00$&$1.00$&$1.00$&$1.00$&$0.98$&$0.44$&$0.19$
\cr
$0.$&{\it MeanMLE}&$100$&$1.01$&$1.01$&$1.01$&$1.00$&$1.01$&$1.01$&$1.01$
\cr
$0.$&{\tty arima}&$100$&$1.01$&$1.01$&$1.01$&$1.01$&$1.01$&$1.01$&$1.01$
\cr\noalign{\smallskip}
$0.$&{\it BLUE}&$200$&$1.00$&$1.00$&$1.00$&$1.00$&$0.99$&$0.58$&$0.26$
\cr
$0.$&{\it MeanMLE}&$200$&$1.00$&$1.00$&$1.00$&$1.00$&$1.00$&$1.00$&$1.00$
\cr
$0.$&{\tty arima}&$200$&$1.00$&$1.00$&$1.00$&$1.00$&$1.00$&$1.00$&$1.00$
\cr\noalign{\smallskip}
$0.5$&{\it BLUE}&$50$&$0.97$&$0.97$&$0.98$&$0.99$&$1.00$&$0.67$&$0.44$
\cr
$0.5$&{\it MeanMLE}&$50$&$0.94$&$0.94$&$0.96$&$1.00$&$1.03$&$0.71$&$0.50$
\cr
$0.5$&{\tty arima}&$50$&$0.94$&$0.94$&$0.96$&$1.00$&$1.03$&$0.74$&$0.55$
\cr\noalign{\smallskip}
$0.5$&{\it BLUE}&$100$&$0.99$&$0.99$&$0.99$&$0.99$&$1.00$&$0.75$&$0.44$
\cr
$0.5$&{\it MeanMLE}&$100$&$0.96$&$0.96$&$0.97$&$0.99$&$1.02$&$0.72$&$0.43$
\cr
$0.5$&{\tty arima}&$100$&$0.96$&$0.96$&$0.97$&$0.99$&$1.02$&$0.74$&$0.44$
\cr\noalign{\smallskip}
$0.5$&{\it BLUE}&$200$&$0.99$&$0.99$&$0.99$&$1.00$&$1.00$&$0.84$&$0.54$
\cr
$0.5$&{\it MeanMLE}&$200$&$0.98$&$0.98$&$0.99$&$1.00$&$1.01$&$0.82$&$0.55$
\cr
$0.5$&{\tty arima}&$200$&$0.98$&$0.98$&$0.99$&$1.00$&$1.01$&$0.84$&$0.56$
\cr\noalign{\smallskip}
$0.9$&{\it BLUE}&$50$&$0.89$&$0.89$&$0.90$&$0.91$&$0.93$&$1.00$&$0.97$
\cr
$0.9$&{\it MeanMLE}&$50$&$0.81$&$0.81$&$0.82$&$0.94$&$0.91$&$1.05$&$1.03$
\cr
$0.9$&{\tty arima}&$50$&$0.80$&$0.80$&$0.81$&$1.58$&$0.91$&$1.08$&$1.02$
\cr\noalign{\smallskip}
$0.9$&{\it BLUE}&$100$&$0.93$&$0.93$&$0.93$&$0.94$&$0.95$&$1.00$&$0.96$
\cr
$0.9$&{\it MeanMLE}&$100$&$0.85$&$0.85$&$0.86$&$0.93$&$0.93$&$1.10$&$1.06$
\cr
$0.9$&{\tty arima}&$100$&$0.85$&$0.85$&$0.86$&$0.93$&$0.93$&$1.09$&$1.04$
\cr\noalign{\smallskip}
$0.9$&{\it BLUE}&$200$&$0.96$&$0.96$&$0.96$&$0.96$&$0.97$&$1.00$&$0.96$
\cr
$0.9$&{\it MeanMLE}&$200$&$0.92$&$0.92$&$0.92$&$0.95$&$0.96$&$1.06$&$0.96$
\cr
$0.9$&{\tty arima}&$200$&$0.92$&$0.92$&$0.92$&$0.95$&$0.96$&$1.04$&$0.97$
\cr\noalign{\smallskip}
$0.95$&{\it BLUE}&$50$&$0.88$&$0.88$&$0.88$&$0.89$&$0.91$&$0.99$&$1.00$
\cr
$0.95$&{\it MeanMLE}&$50$&$0.78$&$0.78$&$0.79$&$0.91$&$0.86$&$1.02$&$1.03$
\cr
$0.95$&{\tty arima}&$50$&$0.77$&$0.77$&$0.79$&$0.91$&$0.91$&$1.01$&$1.01$
\cr\noalign{\smallskip}
$0.95$&{\it BLUE}&$100$&$0.89$&$0.89$&$0.89$&$0.90$&$0.91$&$0.98$&$1.00$
\cr
$0.95$&{\it MeanMLE}&$100$&$0.80$&$0.80$&$0.81$&$0.89$&$0.87$&$1.08$&$1.18$
\cr
$0.95$&{\tty arima}&$100$&$0.80$&$0.80$&$0.81$&$0.89$&$0.87$&$1.08$&$1.14$
\cr\noalign{\smallskip}
$0.95$&{\it BLUE}&$200$&$0.93$&$0.93$&$0.93$&$0.93$&$0.94$&$0.98$&$1.00$
\cr
$0.95$&{\it MeanMLE}&$200$&$0.87$&$0.87$&$0.87$&$0.92$&$0.91$&$1.11$&$1.32$
\cr
$0.95$&{\tty arima}&$200$&$0.87$&$0.87$&$0.87$&$0.92$&$0.91$&$1.13$&$1.29$
\cr\noalign{\smallskip}
\noalign{\smallskip}
\noalign{\hrule}
\end{tabular}
\end{table}

\vfill\eject
\clearpage
\begin{figure}
\centerline{\includegraphics{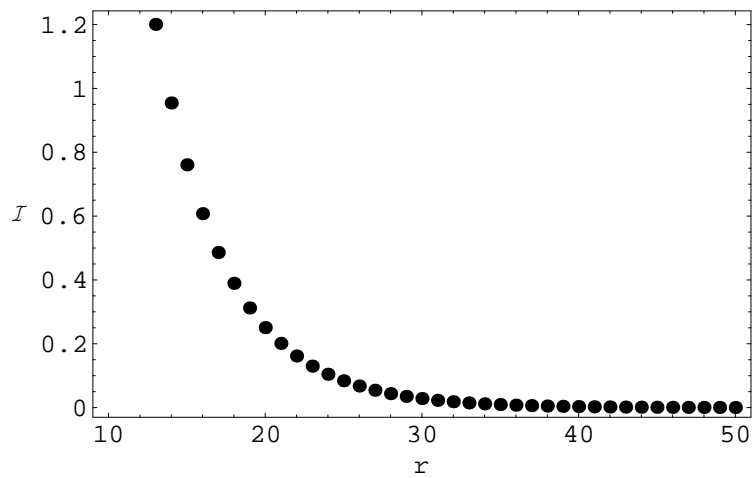}}
\caption{
Kullback-Leibler discrepancy for AR$(r)$ approximation to a $\MA(1)$
with $\theta_1 = 0.9$ and $n=200$.
}
\end{figure}

\vfill\eject
\clearpage
\begin{figure}
\centerline{\includegraphics{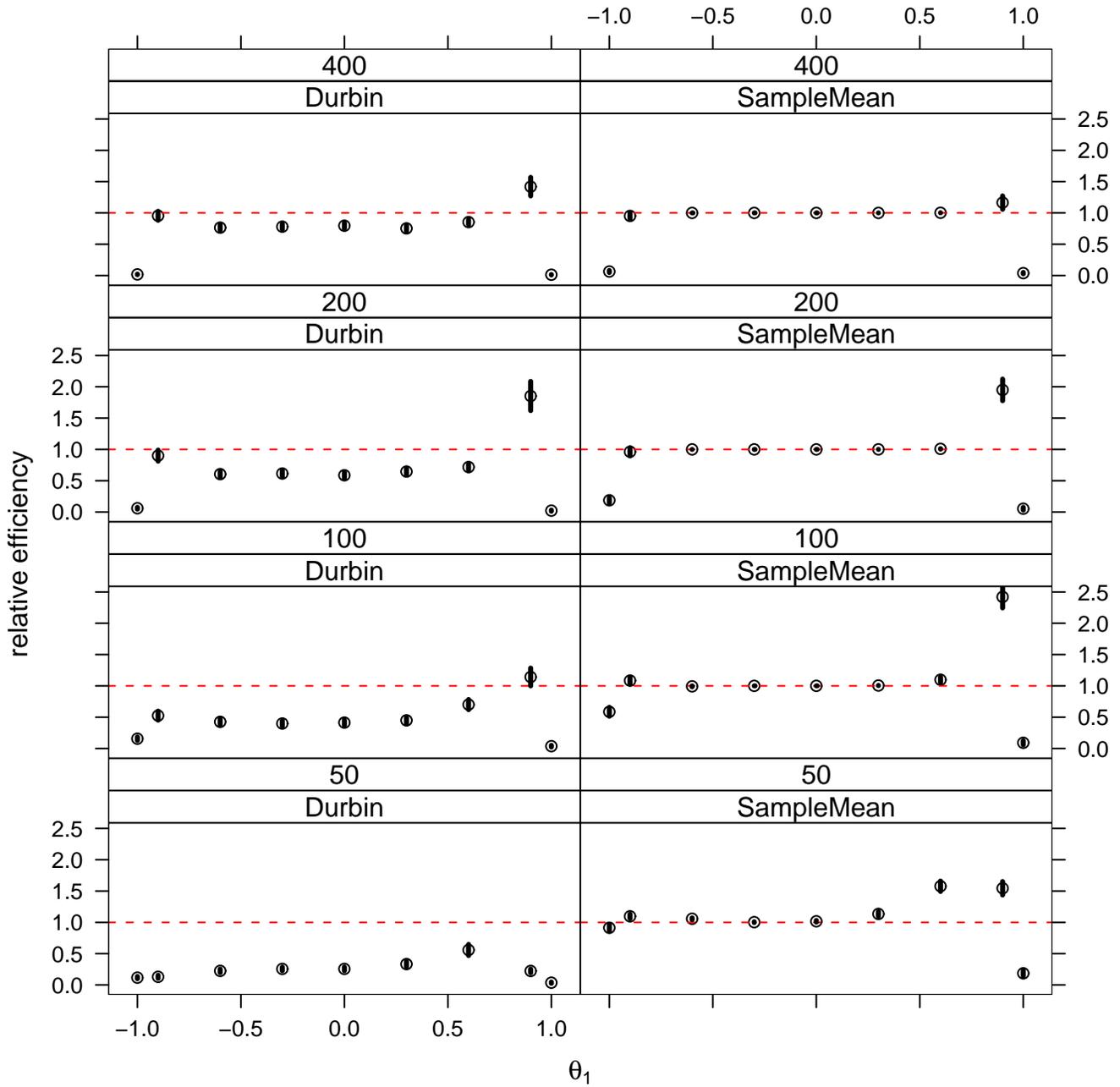}}
\caption{
The vertical lines show the length of the 95\% confidence interval for the statistical efficiency
of {\it SampleMean\/}\  and {\it Durbin\/}\ vs. the MLE based on $10^3$ simulations.
}
\end{figure}

\vfill\eject
\clearpage
\begin{figure}
\centerline{\includegraphics[scale=1.5]{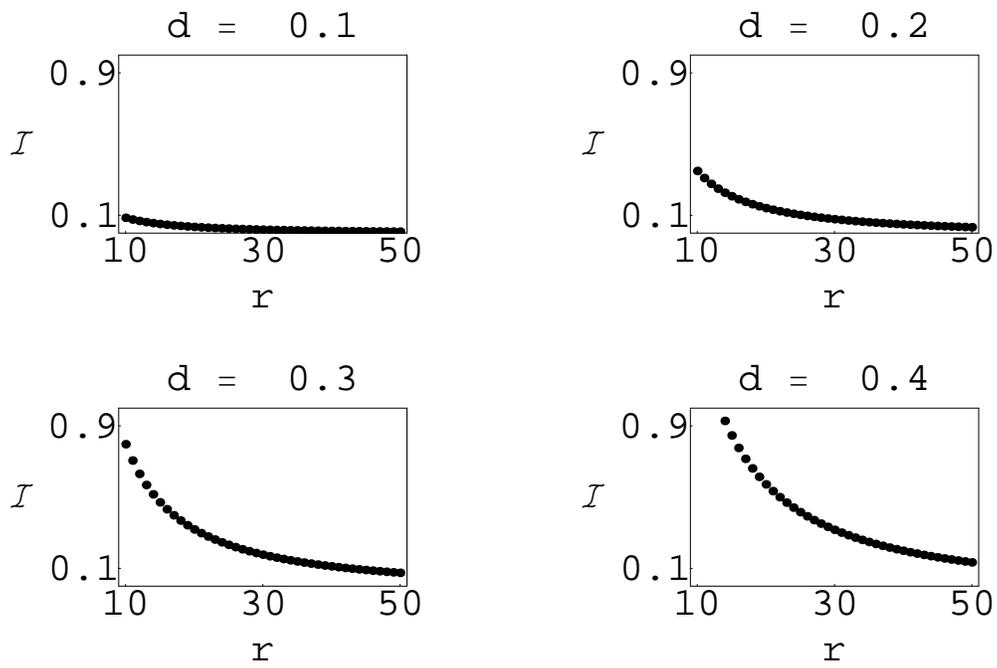}}
\caption{
Kullback-Leibler discrepancy for AR$(r)$ approximation to fractionally
differenced white noise, $\nabla^d z_t = a_t$ for $d=0.1,0.2,0.3,0.4$ and $n=200$.
}
\end{figure}

\end{document}